\documentclass{article}

\usepackage{amsmath}
\usepackage{amsthm}
\usepackage{rotating}

\usepackage[pdftex]{hyperref}

\numberwithin{equation}{section}

\title{Some Combinatorial Identities\\ from the Random Walk}

\author{M.J. Kronenburg}

\begin{document}

\maketitle

\begin{abstract}
Eight combinatorial identities are listed and proved by counting
paths in the one-dimensional random walk.
Four of these identities are assumed to be new.
\end{abstract}

\noindent
\textbf{Keywords}: binomial coefficient, combinatorial identities, random walk.\\
\textbf{MSC 2010}: 05A10, 05A19, 60G50, 68R05

\section{Combinatorial Identities}

Eight combinatorial identities are listed below, of which the last four
are assumed to be new.\\
For nonnegative integer $n$ and complex $m$, $r$:
\begin{equation}\label{combiwalk1}
 \sum_{k=0}^{n}\binom{m+k}{k}\binom{n+r-k}{n-k} = \binom{n+m+r+1}{n}
\end{equation}
\begin{equation}\label{combiwalk2}
 \sum_{k=0}^{n}\frac{m+1}{m+k+1}\binom{m+2k}{k}\binom{2n+r-2k}{n-k} = \binom{2n+m+r+1}{n}
\end{equation}
\begin{equation}\label{combiwalk3}
 \sum_{k=0}^{n}\frac{r+1}{n+r-k+1}\binom{m+2k}{k}\binom{2n+r-2k}{n-k} = \binom{2n+m+r+1}{n}
\end{equation}
\begin{equation}\label{combiwalk4}
\begin{split}
 & \sum_{k=0}^{n}\frac{(m+1)(r+1)}{(m+k+1)(n+r-k+1)}\binom{m+2k}{k}\binom{2n+r-2k}{n-k} \\
={} & \frac{m+r+2}{n+m+r+2}\binom{2n+m+r+1}{n}
\end{split}
\end{equation}
\begin{equation}\label{combiwalk5}
 \sum_{k=0}^{n}\frac{(2k+1)^2}{(n+k+1)(m+k+1)}\binom{2m}{m+k}\binom{2n}{n+k} 
= \frac{1}{n+m+1}\binom{2n+2m}{n+m}
\end{equation}
\begin{equation}\label{combiwalk6}
 \sum_{k=0}^{n}\frac{(k+1)(k+2)}{m+k+2}\binom{2m+k+1}{m}\binom{2n-k}{n} 
= \frac{n+1}{n+m+2}\binom{2n+2m+2}{n+m+1}
\end{equation}
\begin{equation}\label{combiwalk7}
\begin{split}
 & \sum_{k=0}^{n}\frac{(m+1)(m+5)+3(m+2k+1)^2}{(m+k+1)(m+k+2)(m+k+3)}\binom{m+2k}{k} \\
={} &  \frac{4}{n+m+3}\binom{2n+m+2}{n}
\end{split}
\end{equation} 
\begin{equation}\label{combiwalk8}
\begin{split}
 & \sum_{k=0}^{n}\frac{(m+1)(m+5)+(m+2k-1)(m+2k+1)}{(m+k+1)(m+k+2)(m+k+3)}2^{-k}\binom{m+2k}{k} \\
={} & \frac{2^{1-n}}{n+m+3}\binom{2n+m+2}{n}
\end{split}
\end{equation} 
These identities are proved below for integer $m$ and $r$ by counting paths in the
one-dimensional random walk.
Computations show that these identities are also valid for non-integer $m$ and $r$.
The identities (\ref{combiwalk1}) to (\ref{combiwalk4})
are equivalent to two special cases of two identities already
listed and proved in literature \cite{GKP94}.
Changing all $k$ into $n-k$ in the summation term and interchanging $m$ and $r$
transforms (\ref{combiwalk2}) into (\ref{combiwalk3}) and vice versa,
and transforms (\ref{combiwalk1}) and (\ref{combiwalk4}) into itself.
The identities (\ref{combiwalk7}) and (\ref{combiwalk8}) can also be proved by induction on $n$.
The identities can have terms with zeros in numerator and denominator,
in which case the zeros may cancel pairwise by $0/0=1$.
For example taking $m=-1$ in (\ref{combiwalk7}) and (\ref{combiwalk8}) yields:
\begin{equation}
 \sum_{k=0}^{n}\frac{3}{k+2}\binom{2k}{k+1} = \frac{2}{n+2}\binom{2n+1}{n} - 1
\end{equation}
\begin{equation}
 \sum_{k=0}^{n}\frac{4}{k+3}2^{-k}\binom{2k+1}{k+2} = \frac{2^{-n}}{n+3}\binom{2n+4}{n+2} - 2
\end{equation}

\section{Proof of the Combinatorial Identities}

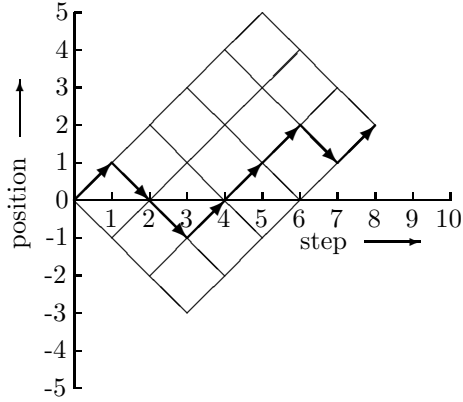
\begin{figure}
\setlength{\unitlength}{0.5cm}
\begin{center}
\begin{picture}(12,11.5)(0.0,0.5)
\put(2,1){\line(0,1){10}}
\put(2,6){\line(1,0){10}}
\multiput(2,1)(0,1){11}{\line(1,0){0.2}}
\multiput(3,6)(1,0){10}{\line(0,1){0.2}}
\put(1.2,0.8){-5}\put(1.2,1.8){-4}\put(1.2,2.8){-3}\put(1.2,3.8){-2}\put(1.25,4.8){-1}
\put(1.5,5.8){0}\put(1.5,6.8){1}\put(1.5,7.8){2}\put(1.5,8.8){3}\put(1.5,9.8){4}\put(1.5,10.8){5}
\put(2.8,5.3){1}\put(3.8,5.3){2}\put(4.8,5.3){3}\put(5.8,5.3){4}\put(6.8,5.3){5}
\put(7.8,5.3){6}\put(8.8,5.3){7}\put(9.8,5.3){8}\put(10.8,5.3){9}\put(11.6,5.3){10}
\put(8,4.7){step}
\put(9.7,4.9){\vector(1,0){1.5}}
\put(0.7,4.7){\begin{rotate}{90}{position}\end{rotate}}
\put(0.55,7.7){\vector(0,1){1.5}}
\multiput(2,6)(1,-1){4}{\line(1,1){5}}
\multiput(2,6)(1,1){6}{\line(1,-1){3}}
\thicklines
\put(2,6){\vector(1,1){1}}\multiput(3,7)(1,-1){2}{\vector(1,-1){1}}
\multiput(5,5)(1,1){3}{\vector(1,1){1}}
\put(8,8){\vector(1,-1){1}}\put(9,7){\vector(1,1){1}}
\thinlines
\end{picture}
\end{center}
\caption{A path in a random walk with $N=8$ and $m=2$.}
\label{figwalk}
\end{figure}

In a one-dimensional random walk, let $N$ be the total number of steps
and let $0$ be the starting position and $m$ be the ending position in steps to the left.
\mbox{$N+m$} and \mbox{$N-m$} are always even, and $m$ can also be negative.
The number of steps to the left was \mbox{$(N+m)/2$} and to the right \mbox{$(N-m)/2$},
because their sum must be $N$ and their difference $m$.
The positions as function of the number of steps of the random walk are therefore always over a
rectangular grid, rotated by $\pi/4$ radians, of $(N+m)/2$ by $(N-m)/2$ (see figure \ref{figwalk}),
and the number of paths over this grid can be summed in different ways,
leading to the combinatorial identities.\\
Let $P(N,m)$ be the total number of paths in this random walk of $N$ steps to position $m$.
This is the number of ways in which \mbox{$(N+m)/2$} steps to the
left can be assigned to the total of $N$ steps:
\begin{equation}\label{walk1}
 P(N,m) = \binom{N}{(N+m)/2}
\end{equation}
In the rectangular grid every path steps from a horizontal or vertical line to the next line once.
Therefore the total number of paths is also a sum over $k$ of paths passing position $r-k$
at step $r+k$ and position $r-k+1$ at step $r+k+1$:
\begin{equation}
 P(N,m) = \sum_{k=0}^{(N-m)/2}P(r+k,r-k)P(N-r-k-1,m-r+k-1)
\end{equation}
or a sum over $k$ of paths passing position $-r+k$ at step $r+k$ and position
$-r+k-1$ at step $r+k+1$:
\begin{equation}
 P(N,m) = \sum_{k=0}^{(N+m)/2}P(r+k,-r+k)P(N-r-k-1,m+r-k+1)
\end{equation}
After substituting (\ref{walk1})
and a change of variables $N$, $m$ and $r$, these equations both result in (\ref{combiwalk1}).\\

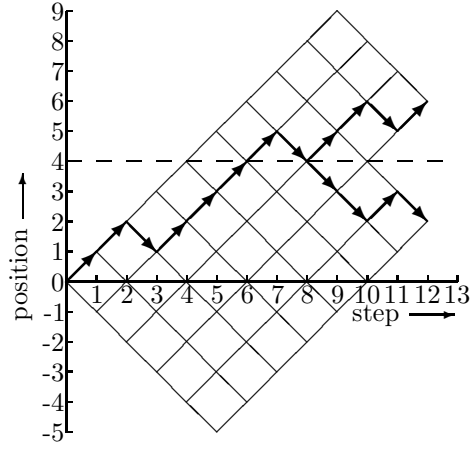
\begin{figure}
\setlength{\unitlength}{0.4cm}
\begin{center}
\begin{picture}(14,16.5)(0.0,0.5)
\put(2,1){\line(0,1){14}}
\put(2,6){\line(1,0){13}}
\multiput(2,1)(0,1){15}{\line(1,0){0.2}}
\multiput(3,6)(1,0){13}{\line(0,1){0.2}}
\put(1.2,0.7){-5}\put(1.2,1.7){-4}\put(1.2,2.7){-3}\put(1.2,3.7){-2}\put(1.25,4.7){-1}
\put(1.5,5.7){0}\put(1.5,6.7){1}\put(1.5,7.7){2}\put(1.5,8.7){3}\put(1.5,9.7){4}\put(1.5,10.7){5}
\put(1.5,11.7){6}\put(1.5,12.7){7}\put(1.5,13.7){8}\put(1.5,14.7){9}
\put(2.8,5.3){1}\put(3.8,5.3){2}\put(4.8,5.3){3}\put(5.8,5.3){4}\put(6.8,5.3){5}
\put(7.8,5.3){6}\put(8.8,5.3){7}\put(9.8,5.3){8}\put(10.8,5.3){9}\put(11.55,5.3){10}
\put(12.55,5.3){11}\put(13.55,5.3){12}\put(14.55,5.3){13}
\put(11.5,4.7){step}
\put(13.4,4.9){\vector(1,0){1.5}}
\put(0.7,4.7){\begin{rotate}{90}{position}\end{rotate}}
\put(0.55,8.1){\vector(0,1){1.5}}
\multiput(2,6)(1,-1){4}{\line(1,1){9}}
\multiput(6,2)(1,-1){2}{\line(1,1){7}}
\multiput(2,6)(1,1){8}{\line(1,-1){5}}
\multiput(10,14)(1,1){2}{\line(1,-1){3}}
\thicklines
\multiput(2,6)(1,1){2}{\vector(1,1){1}}\put(4,8){\vector(1,-1){1}}
\multiput(5,7)(1,1){4}{\vector(1,1){1}}
\put(9,11){\vector(1,-1){1}}\multiput(10,10)(1,1){2}{\vector(1,1){1}}
\put(12,12){\vector(1,-1){1}}\put(13,11){\vector(1,1){1}}
\multiput(10,10)(1,-1){2}{\vector(1,-1){1}}
\put(12,8){\vector(1,1){1}}\put(13,9){\vector(1,-1){1}}
\thinlines
\multiput(2,10)(1,0){13}{\line(1,0){0.5}}
\end{picture}
\end{center}
\caption{$S(N,m,r)=P(N,2r-m)$ with $N=12$, $m=2$, $r=4$.}
\label{figwalk2}
\end{figure}

Let $S(N,m,r)$ be the number of paths in a random walk with nonnegative $m$ passing position $r$
at least once.
Because all paths go from position $0$ to position $m$,
this is equal to (\ref{walk1}) when $0\leq r\leq m$.
When $r\geq m$, each path of $N$ steps to position $m$
and passing position $r$ at least once, corresponds to a path of $N$ steps
to position $2r-m$, where the paths to the last point on position
$r$ are identical, but the paths from there to the ending positions
are mirrored in the line of position $r$ (see figure \ref{figwalk2}).
From this symmetry argument it follows that when $r\leq 0$ or $r\geq m$:
\begin{equation}\label{walk2}
 S(N,m,r) = P(N,2r-m) = \binom{N}{(N+m)/2-r}
\end{equation}
When \mbox{$r=0$} or \mbox{$r=m$}, the positions of the starting or ending point,
this becomes equal to (\ref{walk1}), and when \mbox{$r=(N+m)/2$}
or \mbox{$r=-(N-m)/2$}, the extreme positions, this becomes $1$.
For any $r$ the following recursion relation is fulfilled:
\begin{equation}\label{walkrecur}
 S(N+1,m+1,r) = S(N,m,r) + S(N,m+2,r)
\end{equation}
which in words means that for any $r$ the number of paths to a position
is the sum of the number of paths of one step earlier to the position
one to the left and one to the right.\\

Let $T(N,m,r)$ be the number of paths in a random walk with nonnegative $m$
without passing position $r$, where $r$ is negative:
\begin{equation}
 T(N,m,r) = P(N,m) - S(N,m,r)
\end{equation}
Then it follows with (\ref{walk1}) and (\ref{walk2}) that for example,
with $\mu=(N+m)/2$:
\begin{equation}\label{walk3a}
 T(N,m,-1) = \dfrac{m+1}{\mu+1}\binom{N}{\mu}
\end{equation}
\begin{equation}\label{walk3b}
 T(N,m,-2) = \dfrac{(m+2)(N+1)}{(\mu+1)(\mu+2)}\binom{N}{\mu}
\end{equation}
\begin{equation}\label{walk3c}
 T(N,m,-3) = \dfrac{(m+3)[(m+1)(m+5)+3(N+1)^2]}{4(\mu+1)(\mu+2)(\mu+3)}\binom{N}{\mu}
\end{equation}
\begin{equation}\label{walk3d}
 T(N,m,-4) = \dfrac{(m+4)(N+1)[(m+2)(m+6)+N(N+2)]}{2(\mu+1)(\mu+2)(\mu+3)(\mu+4)}\binom{N}{\mu}
\end{equation}
In the rectangular grid every path reaches position $r$, where $0\leq r\leq m$,
at some number of steps for the first time.
Therefore the total number of paths is also
a sum over $k$ of paths reaching position $r$ for the first time after \mbox{$r+2k$} steps,
which is the product of the number of paths to position $r-1$ after $r+2k-1$ steps without
passing position $r$, and the number of paths from position $r$ to position $m$
in $N-r-2k$ steps:
\begin{equation}\label{walk4}
 P(N,m) = \sum_{k=0}^{(N-m)/2}T(r+2k-1,r-1,-1)P(N-r-2k,m-r)
\end{equation}
or a sum over $k$ of paths passing that position for the last time after $r+2k$ steps,
which is the product of the number of paths to position $r$ after $r+2k$ steps,
and the number of paths from position $r+1$ to position $m$ in $N-r-2k-1$ steps
without passing position $r$:
\begin{equation}\label{walk5}
 P(N,m) = \sum_{k=0}^{(N-m)/2}P(r+2k,r)T(N-r-2k-1,m-r-1,-1)
\end{equation}
After substituting (\ref{walk1}) and (\ref{walk3a}),
and a change of variables $N$, $m$ and $r$, these equations result in
(\ref{combiwalk2}) and (\ref{combiwalk3}).
When $r<0$ or $r>m$ similar equations apply, leading to the same result.
When $r<0$ and in addition position $r-1$ is not passed, the following equation applies:
\begin{equation}\label{walk6}
\begin{split}
 & S(N,m,r) - S(n,m,r-1) \\
={} & \sum_{k=0}^{(N-m)/2+r}T(2k-r-1,-r-1,-1)T(N-2k+r,m-r,-1)
\end{split}
\end{equation}
which results in (\ref{combiwalk4}).
The number of paths not passing position $-1$ is also a sum
of paths passing position $2k$ at step number $r$, where $r$ is even:
\begin{equation}
 T(N,0,-1) = \sum_{k=0}^{r/2}T(r,2k,-1)T(N-r,2k,-1)
\end{equation}
which results in (\ref{combiwalk5}),
or a sum of paths passing position $r-k$ at step $r+k$ and position
$r-k+1$ at step $r+k+1$:
\begin{equation}
 T(N,0,-1) = \sum_{k=0}^{r}T(r+k,r-k,-1)T(N-r-k-1,r-k+1,-1)
\end{equation}
which results in (\ref{combiwalk6}),
or a sum of paths reaching position $2$ for the first time at step $2k+2$
and not passing position $-1$:
\begin{equation}
 T(N,m,-1) = \sum_{k=0}^{(N-m)/2}T(N-2k-2,m-2,-3)
\end{equation}
which results in (\ref{combiwalk7}),
or a sum of paths reaching position $3$ for the first time at step $2k+3$
and not passing position $-1$:
\begin{equation}
 T(N,m,-1) = \sum_{k=0}^{(N-m)/2} T(N-2k-3,m-3,-4) 2^k
\end{equation}
which results in (\ref{combiwalk8}).

\pdfbookmark[0]{References}{}

\end{document}